# Transactive Control of Air Conditioning Loads for Mitigating Microgrid Tie-line Power Fluctuations

Yao Yao, Peichao Zhang
Department of Electrical Engineering
Shanghai Jiao Tong University
Shanghai, China

*Abstract*—This paper presents a distributed control strategy for air conditioning loads (ACLs) to participate in the scheme of mitigating microgrid tie-line power fluctuations. The concept of baseline load is emphasized for ACL control in this paper. To obtain the target aggregated power of ACLs, an algorithm based on the principle of low-pass filter (LPF) is derived. For better robustness of the control strategy, feedback of the states of indoor temperatures is introduced for baseline load correction. A transactive control method is then put forward to allocate the target aggregated power to each ACL. This method can satisfy customers' differentiated requirements for comfort, protect the privacy and enhance the security of the controlled appliances. For the microgrid control center, it can simplify the downlink control and avoid measuring the power of uncontrolled loads which reduces the implementation cost. Simulation results shows ACLs can effectively provide microgrid tie-line smoothing services using the proposed method.

## I. INTRODUCTION

The integration of large scale intermittent sources, such as wind power and photovoltaic power, has adverse effect on power quality and power system stability. With the increase of the penetration of renewable energy, it has been an important research subject how to mitigate the power fluctuation. The service of smoothing the fluctuations has been traditionally provided by energy storage devices. However, the current storage technologies have high cost [1].

Many recent papers have focused on the use of thermostatically controlled loads (TCLs) such as air-conditioners (AC), heat pumps and water heaters for their ability of thermal energy storage [1]-[5]. They are demand response resources of great potential because they can be transformed to virtual storages of large amount, low cost and fast response speed through certain control means.

The key point of smoothing fluctuations is realizing active power tracking. So far, there have been many effective control strategies with TCLs. In [2], a state-queueing model [3] is adopted and target power tracking is realized by direct on-off control of air conditioning loads (ACLs). Model predictive control (MPC) and a control strategy with setpoint adjustment is proposed in [4]. In [5], a continuum-scale partial differential equation model for aggregate TCLs is derived, based on which a sliding mode controller is developed to smooth the fluctuations.

However, the existing methods have several problems:

1) The control algorithms need to measure or predicate the uncontrollable loads [2], [4]. The loads' characteristics of volume, variety and disperstiveness make it costly to implement.

2) Many control strategies use the direct load control, including switching control [2] and setpoint adjustment [4], [5]. Although faster response can be obtained, requirements for communication bandwidth can be high when the scale is large. What's more, customers will face more security risks when the direct switch control of appliances is used.

3) Some control strategies [2]-[5] need to collect customers' preferences or the thermal parameters / models of the buildings, which is very costly in practice and might cause privacy issues.

In order to improve user experience and reduce implementation cost, a transactive control method for ACLs to smooth microgrid tie-line power fluctuations is put forward in this paper. Transactive control is also called market-based control (MBC). It adopts the market equilibrium mechanism of micro economics. Optimal allocation problem of finite resources in the computer field has been solved by MBC for a long time [6]. In recent years, it has also been used to implement the coordination control of large-scale distributed energy in some demonstration projects in USA and Europe [7].

The control method in this paper consists of two parts. The first part aims at the whole ACL cluster. The control target of the aggregated power of ACLs is derived based on the low-pass filter (LPF) principle. The second part aims at the ACL individual. The allocation of the control target to each ACL is accomplished by transactive control through establishing a virtual market in the microgrid.

The rest of this paper is organized as follows: Section II presents an algorithm for tie-line power smoothing. Section III introduces the allocation strategy of target aggregated power to each ACL. Section IV examines the smoothing effects. The conclusions and future work are summarized in Section V.

## II. TARGET AGGREGATED POWER OF ACLS

### A. Smoothing Strategy of Microgrid Tie-line Power Fluctuations

Fig. 1 is a microgrid system, where $P_w$ is the wind power, $P_{AC}$ is the aggregated power of ACLs participating in the tie-line power smoothing scheme, $P_L$ is the total power of the uncontrollable loads, and $P_g$ is the tie-line power. Ignoring the line loss, we can get the equation below at time $k$:

$$P_g[k] = P_{AC}[k] + P_L[k] - P_W[k]. \quad (1)$$

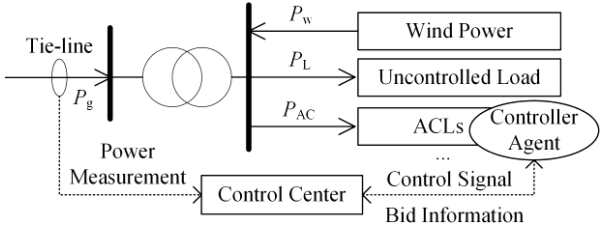

**Fig. 1. Schematic of a microgrid**

The free load of ACLs without external control is called AC baseline load ($P_{ACbase}$) in this paper. According to Fig.1, the free tie-line power with all ACLs uncontrolled is

$$P_{g0}[k] = P_{ACbase}[k] + P_L[k] - P_W[k], \quad (2)$$

the fluctuations in which are mainly caused by renewable energy.

LPF is used to smooth the tie-line power fluctuations. The tie-line target power can be calculated using the recursive form:

$$P_{gLPF}[k] = \alpha P_{gLPF}[k-1] + (1-\alpha) P_{g0}[k] \quad (3)$$

where $\alpha = \tau/(\tau + \Delta t)$ is the filter coefficient, $\tau$ is the time constant, and $\Delta t$ is the control cycle.

To track the tie-line target power, the aggregated power of ACLs needs to make adjustment as below:

$$\Delta P_{AC}[k] = P_{gLPF}[k] - P_{g0}[k]. \quad (4)$$

$\Delta P_{AC}$ contains high frequency fluctuating components in the tie-line power. Thus, the target power of ACLs, $P_{AC}^*$, can be calculated as:

$$P_{AC}^*[k] = P_{ACbase}[k] + \Delta P_{AC}[k]. \quad (5)$$

For ACLs, the baseline load in the absence of demand response events should be satisfied first, based on which the adjustment $\Delta P_{AC}[k]$ is then added. Although the ACL could be treated as a virtual battery [8], the control strategy needs good baseline load estimation to make sure the indoor temperatures are always kept within limits, which significantly differs from the battery storage as the baseline load is normally not a concern for the latter.

### B. Baseline Load Estimation

The methods of AC baseline load estimation have been discussed in [4], [7], [9]. Due to space limitations, this paper will not study this topic. Instead, more efforts will be made to correct the estimation. Thus, this paper simply employs the multiple quadratic regression method to calculate the estimated baseline load, $P_{ACbase0}$. According to the ACL model [10], [11], the main parameters including the outdoor temperature, solar radiation [11], [12] and total rated power of ACLs are considered in the regression model.

### C. AC Baseline Load Correction Based on SOA Feedback

Baseline load estimation error inevitably exists, resulting in deviation from ideal range of the indoor temperature. A correction method using state of indoor temperature of air-conditioner (SOA) as feedback is put forward in this paper.

To quantify the regulation capacity of ACLs and the comfort levels of customers, SOA is defined as:

$$SOA = \begin{cases} \dfrac{T_{air} - T_{set}}{T_{max} - T_{set}}, & T_{air} \geq T_{set} \\ \dfrac{T_{air} - T_{set}}{T_{set} - T_{min}}, & T_{air} < T_{set} \end{cases} \quad (6)$$

where $T_{set}$ is the indoor temperature setpoint, $T_{max}$ and $T_{min}$ is the upper and lower temperature limit, and $T_{air}$ is the current indoor temperature.

It's obvious that $SOA \in [-1,1]$. The closer it is to zero, the more regulation capacity there is. If it is close to 1 or -1, the indoor temperature is near the upper/lower limit.

$S$ is used to measure the general state of indoor temperatures of the ACL cluster, which is defined as:

$$S = \frac{1}{n} \sum_{i=1}^{n} SOA_i \quad (7)$$

where $n$ is the number of ACLs to be controlled.

Taking $S$ as feedback, the correction formula is:

$$\begin{cases} P_{ACbase}[k] = P_{ACbase0}[k] + P_{adj}[k] \\ P_{adj}[k] = \Delta P_{adj}[k]\% \times P_{ACbase0}[k] + P_{adj}[k-1]e^{-\gamma} \end{cases} \quad (8)$$

where $P_{ACbase0}$ is the original estimated baseline load, and $P_{adj}[k]$ is the correction of the $k^{th}$ control cycle, consisting of a proportional component and an attenuation component. $\Delta P_{adj}[k]$ is the proportional coefficient determined by $S$ and $\gamma > 0$ is the attenuation coefficient. When $\Delta P_{adj}[k] \neq 0$, the attenuation component can accelerate the adjustment speed, and when $\Delta P_{adj}[k] = 0$, it will gradually reduce the correction amount.

The control flow of the correction method is shown in Fig. 2.

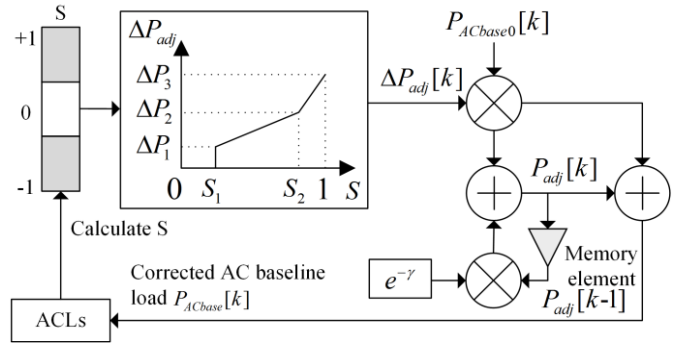

**Fig. 2. AC baseline load correction method based on SOA feedback**

$\Delta P_{adj}(S)$ is an odd function. Only the first quadrant part is shown in Fig. 2 and expressed as below:

$$\Delta P_{adj} = \begin{cases} \Delta P_1 + \dfrac{\Delta P_2 - \Delta P_1}{S_2 - S_1}(S - S_1), S_1 \leq S < S_2 \\ \Delta P_2 + \dfrac{\Delta P_3 - \Delta P_2}{S_3 - S_2}(S - S_2), S_2 \leq S \leq 1 \end{cases}. \quad (9)$$

## III. ALLOCATION OF AGGREGATED POWER OF ACLs

So far, the target aggregated power $P^*_{AC}[k]$ has been obtained. Another key point is how to allocate $P^*_{AC}[k]$ to the distributed ACLs. The following principles are considered when designing the allocation method: (1) satisfying the different requirements for comfort; (2) no parameters or model information of buildings and ACLs are needed; (3) no explosion of switch control rights of appliances; (4) suitable for all kinds of TCLs; (5) supporting plug-and-play control of ACLs.

To achieve the above goals, a virtual market is established in the microgrid, and the allocation is accomplished by transactive control based on the market equilibrium principle.

### A. General Control Flow

The overall control flow chart is shown in Fig. 3.

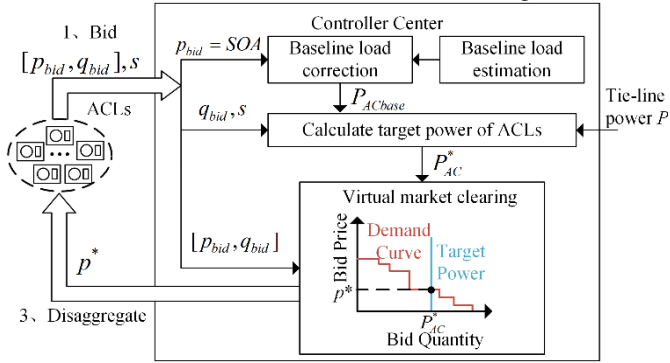

**Fig. 3. Overall control flow chart**

As Fig. 3 shows, each control cycle includes three stages:

1) Bidding stage: Each ACL sends bidding information to the microgrid control center (MGCC) before the next control cycle.

2) Aggregation stage: First, the MGCC calculates the target power $P^*_{AC}[k]$. Then, the virtual market collects and aggregates the bids of ACLs, forming the demand curve. Finally, the MGCC clears the virtual market by solving the intersection of $P^*_{AC}[k]$ and the demand curve.

3) Disaggregation stage: The virtual market broadcasts the clearing result $p^*$ and accomplishes the allocation of $P^*_{AC}[k]$.

### B. Bidding Strategy of ACLs

Each ACL has a controller as its agent. The bidding information of the $i^{th}$ ACL at the $k^{th}$ control cycle is:

$$B_i[k] = ([p_{bid}, q_{bid}], s)_i[k] \quad (10)$$

where bid price $p_{bid}$=SOA which only plays a role as control signal without economic meaning; bid quantity $q_{bid}$ is the power of ACL in operation which is normally set as the rated power [2]-[4], and $s$ denotes the operation state at the time of bidding which equals 1 when the ACL is on and 0 when it's off.

Thermal parameters and comfort preferences of customers are not needed under the bidding mechanism above, which keeps the private information at the user end and protects the privacy.

### C. Virtual Market Clearing

As is shown in Fig. 4, the virtual market sorts the bid information in descending order of bid price $p_{bid}$, forms the demand curve and solves the clearing price $p^*$ at the intersection of demand curve and target power $P^*_{AC}[k]$. In clearing scenario as shown in Fig. 4(b), $p^*=(p_{bid1}+p_{bid2})/2$.

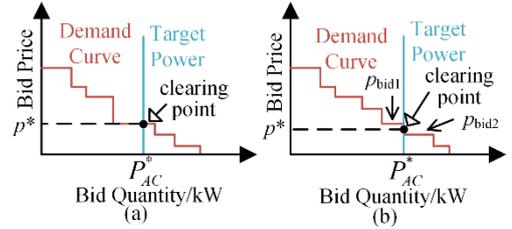

**Fig. 4. Schematic of market clearing**

According to formulas (2)-(5), the value of $P_L$-$P_W$ is needed to calculate the target power $P^*_{AC}[k]$, which can be obtained by:

$$P_L[k] - P_W[k] = P_g[k] - \sum_{i=1}^{n} s_i[k]q_{bid,i}[k] \quad (11)$$

where $s_i$ and $q_{bid,i}$ is the operation state and bid quantity of the $i^{th}$ ACL, and $n$ is the total number of controlled ACLs.

Thus, only the tie-line power $P_g$ is needed to be measured with this method, which can significantly reduce the implementation cost.

### D. Response to The Clearing Result

Clearing price $p^*$ is the only control signal sent to each ACL. To track the target power $P^*_{AC}[k]$, each ACL should be turned off if its bid price is lower than $p^*$, otherwise it should be turned on. A setpoint adjustment method is used:

$$T_{set} = \begin{cases} T_{max} - \varepsilon, p^* \geq p_{bid} \\ T_{min} + \varepsilon, p^* < p_{bid} \end{cases} \quad (12)$$

where $\varepsilon$ is determined by each ACL to make sure the indoor temperature is in the range of $[T_{min}, T_{max}]$ during the control cycle.

Because there is no need for the MGCC to specify the operation state or the setpoint of each ACL, this method effectively simplifies the downlink control. Compared with the direct switch control as adopted in [2], [3], the indirect control by clearing price has lower cyber security risks as the indoor temperatures are finally maintained by the local controller as shown in (12).

## IV. SIMULATION RESULTS

### A. Parameter Settings of The Simulation Case

In the simulation case, there are 450 ACLs to be controlled in the microgrid. The proportion of ACLs at peak load is about 40% and the wind power makes up about 27% (the ratio of the installed capacity to the load peak).

Second-order ETP model [10], [11] is adopted for modeling the ACLs. The simulation step is 5s while the data updating & recording cycles are 10s. The control cycle is 1min, and the ACL controllers bid 5s before the next control cycle.

Main parameter settings are shown in Table I-Table III where U(a,b) denotes the uniform distribution and N(avg,std) denotes the normal distribution.

TABLE I. MAIN PARAMETER SETTINGS OF ACLS

| Area/ $m^2$ | Air Change Freq/(Times/h) | Window-Wall Ratio | SHGC | EER |
|---|---|---|---|---|
| U(88,176) | N(0.5,0.06) | N(0.15,0.01) | U(0.22,0.5) | U(3,4) |
| $R_{th}$ of Roof/ (ºC.m²/W) | $R_{th}$ of Wall/ (ºC.m²/W) | $R_{th}$ of Floor/ (ºC.m²/W) | $R_{th}$ of Window/ (ºC.m²/W) | $R_{th}$ of Door/ (ºC.m²/W) |
| N(5.28,0.70) | N(2.99,0.35) | N(3.35,0.35) | N(0.38,0.03) | N(0.88,0.07) |

Note: SHGC denotes solar heat gain coefficient, EER denotes energy efficiency ratio, and $R_{th}$ denotes thermal resistance.

TABLE II. MAIN PARAMETER SETTINGS OF CONTROLLERS

| Deadband of Thermostat /°C | $T_{set}$/°C | $T_{high}$/°C[1] | $T_{low}$/°C[2] |
|---|---|---|---|
| U(0.2,0.4) | N(26,0.5) | U(2,3) | U(2,3) |

Note 1: $T_{high}= T_{max} - T_{set}$; 2: $T_{low}= T_{set} - T_{min}$.

TABLE III. MAIN PARAMETER SETTINGS OF MGCC

| $\tau$/min | $S_1$ | $S_2$ | $\Delta P_1$/% | $\Delta P_2$/% | $\Delta P_3$/% | $\gamma$ |
|---|---|---|---|---|---|---|
| 50 | 0.5 | 0.8 | 1 | 2 | 3 | 0.02 |

Uncontrollable load, wind power, outdoor temperature and solar radiation [12] on the simulation day are shown in Fig. 5.

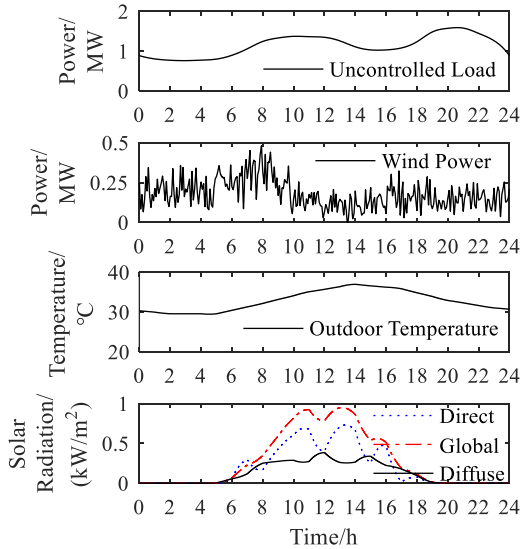

Fig. 5. Uncontrollable load, wind power, outdoor temperature and solar radiation

The MGCC collects the real-time environmental data and calculates the AC baseline load estimated value $P_{ACbase0}$ with the multiple quadratic regression model. The real power of uncontrolled ACLs and the estimated baseline $P_{ACbase0}$ are compared in Fig. 6.

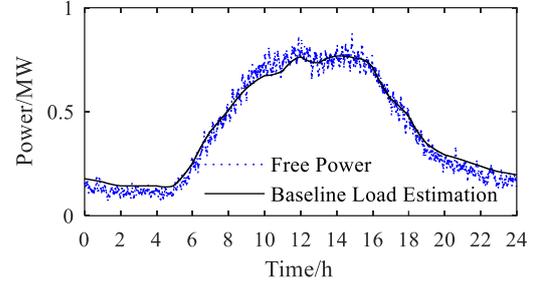

Fig. 6. Free power and baseline load estimation of ACLs

### B. Effect of Tie-line Power Smoothing

The smoothing effect of tie-line power fluctuations with this method is shown in Fig. 7, where $P_{gLPF}$ is the target tie-line power, $P_g$ is the actual tie-line power, and $P_{g0}$ is the free tie-line power when ACLs are not controlled.

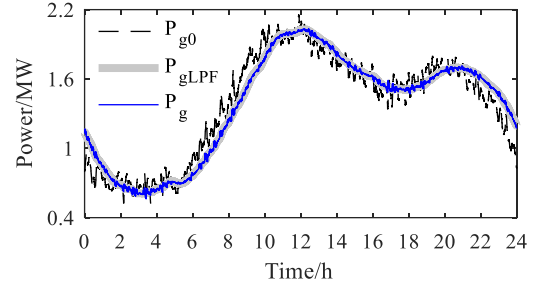

Fig. 7. Smoothing effect

Define 10min fluctuation rate of tie-line power at time $t$ as below:

$$R_{fluc10min}^t = \max_{i \in (t-10/t_r, t]} (P_g^i) - \min_{i \in (t-10/t_r, t]} (P_g^i) \quad (13)$$

where $t_r$ is the record cycle (min), and $P_g^i$ is the actual tie-line power at the $i^{th}$ min.

The 10min fluctuation rates of tie-line power with / without controlling the ACLs are compared in Fig. 8.

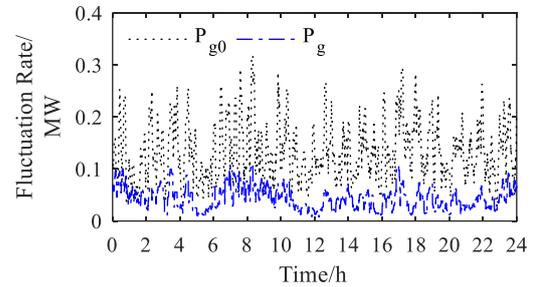

Fig. 8. Fluctuation rate

Fig. 7 and Fig. 8 show that the actual tie-line power tracks the target power well and the fluctuations are effectively smoothed with this control method.

### C. Effect of SOA Feedback Control

To demonstrate the effect of SOA feedback control under circumstance of large estimation errors, the AC baseline load estimated value is intentionally adjusted by ±10% based on the value in Fig. 6.

The simulation results with and without SOA control are shown in Fig. 9 and Fig. 10.

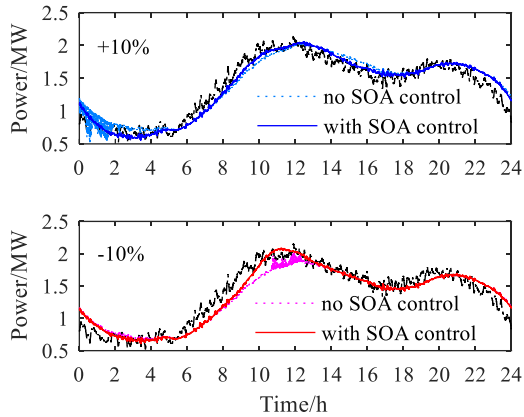

Fig. 9. Smoothing effect with/without SOA feedback control

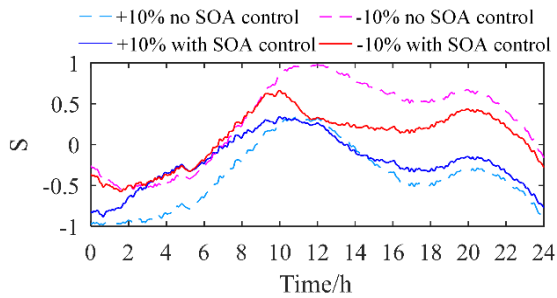

Fig. 10. Changes of $S$ with/without SOA feedback control

As is shown in Fig. 10, when SOA feedback control is not applied, $S$ reaches the upper/lower limit during some periods. Because the comfort of customers is of top priority in this method, the ACLs cannot track the target power accurately and control failure occurs. If the estimated power is too large, it's easy for $S$ to reach the lower limit when the outdoor temperature is low and the solar radiation is weak (e.g. 0:00 a.m to 6:00 a.m); if the estimated power is too small, it's easy for $S$ to reach the upper limit when the outdoor temperature is high and solar radiation is strong (e.g. 9:00 a.m to 2:00 p.m). With SOA feedback control, the estimated value during these periods is corrected according to $S$, which ensures that $S$ is within the ideal range.

## V. CONCLUSION AND PROSPECT

To reduce the size of the expensive battery storage, this paper proposes a distribution control method for ACLs to provide the tie-line power smoothing service. Simulation results show that the power fluctuations can be significantly mitigated with this control strategy.

In contrast with the control of battery storages, baseline load estimation is an important subject for the control of ACLs. The baseline load correction method based on SOA feedback can reduce the requirements for estimation accuracy and effectively enhance the robustness of the control strategy.

The control strategy based on the idea of transactive control can provide high scalability, simplify the downlink control and satisfy customers' different requirements for comfort while protecting the privacy and enhancing the security of the controlled appliances.

When realizing the control method, only the single power of the microgrid tie-line needs to be measured, which effectively reduces the deployment and operation costs of the proposed method.

Future work will focus on the coordination control of TCLs and battery storages to improve the control effect.